\documentclass[a4paper,12pt]{article}
\usepackage{amssymb,amsmath}
\newcommand{\dz}{\delta_0}
\newcommand{\dw}{\delta_1}
\renewcommand{\b}[1]{\mathbf{#1}}
\newcommand{\cl}[1]{\mathcal{#1}}

\newcommand{\Z}{\mathbb{Z}}
\newcommand{\N}{\mathbb{N}}

\newcommand{\twosum}[2]{\sum_{\substack{#1\\#2}}}
\newcommand{\ep}{\varepsilon}
\newcommand{\beql}[1]{\begin{equation}\label{#1}}
\newcommand{\eeq}{\end{equation}}
\newcommand{\meas}{\text{meas}}
\newcommand{\andd}{\;\;\;\mbox{and}\;\;\;}

\newtheorem{theorem}{Theorem}

\newtheorem{lemma}{Lemma}

\begin{document}

\title{The Differences Between Consecutive Smooth Numbers}

\author{D.R. Heath-Brown\\Mathematical Institute, Oxford}

\date{}

\maketitle

\begin{flushright}
  {\em In Celebration of the Seventy-Fifth Birthday of Robert Tijdeman}
\end{flushright}


\section{Introduction}

If $y$ is a positive real number, an integer $n\in\N$ is said to be
$y$-smooth if all its prime factors $p\mid n$ satisfy $p\le y$.  We
are interested in upper bounds for large gaps between consecutive
$y$-smooth numbers.  For example, if $a_1,a_2,\ldots$ is the sequence
of $y$-smooth numbers in increasing order, one might ask for upper
bounds for
\[\max\{a_{n+1}-a_n:\; a_n\le x\}\]
in terms of $x$ and $y$. In this paper we will primarily be interested
in a measure for the frequency of large gaps, given by
\[\sum_{a_n\le x}(a_{n+1}-a_n)^2.\]
The reader should note that there is a dependence on $y$ which is not
mentioned explicitly above. In order to assess bounds for this latter
sum we will want to know how many $y$-free integers $a_n\le x$ there
are. The notation $\psi(x,y)$ is standard for this quantity.  It has
been extensively investigated, but for our purposes it will suffice to
know that
\[ x\ll_{\ep}\psi(x,x^{\ep})\le x\]
for any fixed $\ep>0$, and that
\[\psi(x,y)=x^{1+o(1)}\]
if $y\ge (\log x)^{f(x)}$ for some function $f(x)$ tending to infinity
with $x$ (see Hildebrand and Tenenbaum \cite[Corollary 1.3]{HT}, for
example). When $\psi(x,y)=x^{1+o(1)}$ one might hope that
\beql{conj}
\sum_{a_n\le x}(a_{n+1}-a_n)^2\ll_{\ep} x^{1+\ep}
\eeq
for any fixed $\ep>0$.

If $y$ is not too small compared to $x$ it turns out that questions
about $y$-smooth numbers are no harder, and sometimes easier, than the
corresponding questions about primes.  For example, recent work of
Matom\"{a}ki and Radziwi\l\l\, \cite[Corollary 1]{MR} shows that for any
$\ep>0$ the gaps between consecutive $x^{\ep}$-smooth numbers up to
$x$ are at most $O_{\ep}(x^{1/2})$.  The corresponding result for
primes is just out of reach, even on the Riemann hypothesis. 

There has
been much work on the sum
\[\sum_{p_n\le x}(p_{n+1}-p_n)^2.\]
Under the Riemann hypothesis, work of Selberg \cite{Selb} shows that
\[\sum_{p_n\le x}(p_{n+1}-p_n)^2\ll x(\log x)^3,\]
while Yu \cite {GY} obtains
\[\sum_{p_n\le x}(p_{n+1}-p_n)^2\ll_{\ep} x^{1+\ep},\]
subject only to the Lindel\"{o}f hypothesis. An unconditional bound of
the same strength seems far out of reach, and the present paper is
therefore designed to investigate the extent to which one can
establish the corresponding bound (\ref{conj}) unconditionally.
We shall content ourselves with an investigation of $x^{\ep}$-smooth
numbers.  However our approach generalizes to $y$-smooth number in
general, so long as $(\log x)/(\log y)$ is at most a small power of
$log x$.

\begin{theorem}\label{t1}
Let $a_n$ be the $x^{\ep}$-smooth numbers, in increasing order.  Then
\[\twosum{a_n\le x}{a_{n+1}-a_n\ge x^{1/3+\ep}}(a_{n+1}-a_n)^2\ll_{\ep}
  x^{1+\ep}.\]
\end{theorem}

Thus we achieve (\ref{conj}) for those gaps that have length at least
$x^{1/3+\ep}$. Unfortunately however our new method breaks down entirely
for smaller gaps.  To cover the remaining range we can use
pre-existing methods, which lead to our next result.

\begin{theorem}\label{t2}
Let $a_n$ be the $x^{\ep}$-smooth numbers, in increasing order.  Then
\[\sum_{a_n\le x}(a_{n+1}-a_n)^{3/2}\ll_{\ep} x^{1+\ep}.\]
\end{theorem}

As a corollary we then obtain the following bound.

\begin{theorem}\label{t3}
Let $a_n$ be the $x^{\ep}$-smooth numbers, in increasing order.  Then
\[\sum_{a_n\le x}(a_{n+1}-a_n)^2\ll_{\ep} x^{7/6+\ep}.\]
\end{theorem}

The key idea behind our proof of Theorem \ref{t1} 
is a new type of estimate for a certain
mean value of Dirichlet polynomials.  Let $T\ge 20$
and let $\cl{M}$ be a set of
distinct integers $m\in(0,T]$.   We will write $R=\#\cl{M}$. For each
  $m\in\cl{M}$ let $\ep_m$ be a complex number of modulus at most 1.
Suppose that $N$ is a positive integer
and that $q_1,\ldots,q_N$ are real coefficients in $[0,1]$.  Write
\[M(s):=\sum_{m\in\cl{M}}\ep_m m^{-s}\andd Q(s):=\sum_{n\le N}q_n n^{-s}.\]
The mean value in which we are interested is then
\[\cl{I}(\cl{M},Q):=\int_0^T|M(it)Q(it)|^2dt,\]
which will be related to
\[\cl{J}(\cl{M},Q):=\twosum{(m_1,m_2,n_1,n_2)\in\cl{M}^2\times\N^2}
     {|\log(m_1n_1/m_2n_2)|\le 2\pi/T}q_{n_1}q_{n_2}.\]
We then have the following results.
\begin{theorem}\label{key}

\begin{enumerate}
\item[(i)]
  We have
  \[\cl{I}(\cl{M},Q)\ll T\cl{J}(\cl{M},Q).\]
  Moreover, if $\ep_m=1$ for every $m\in\cl{M}$ then
  \[\cl{I}(\cl{M},Q)\ll T\cl{J}(\cl{M},Q)\ll \cl{I}(\cl{M},Q).\]
\item[(ii)]
Under the Lindel\"{o}f Hypothesis, for any $\eta>0$ and any $Q(s)$ we have
\beql{LHB}
\cl{I}(\cl{M},Q)\ll_{\eta}N^2R^2+(NT)^{\eta}NRT.
\eeq
\item[(iii)]
Unconditionally, for any $\eta>0$ and any $Q(s)$ we have
\[\cl{I}(\cl{M},Q)\ll_{\eta}N^2R^2+(NT)^{\eta}\{NRT+NR^{7/4}T^{3/4}\}.\]
Moreover, (\ref{LHB}) holds when either $N\ge T^{2/3}$ or $R\le T^{1/3}$.
\item[(iv)]
  If
  \[Q(s)=\frac{1}{k!}\left(\sum_{N^{1/k}/2<p\le N^{1/k}}p^{-s}\right)^k,\]
  with $p$ running over primes, then
  \[\cl{I}(\cl{M},Q)\ll_k N^2R^2+ R^2T+NRT.\]
\end{enumerate}
\end{theorem}

Part (ii) is included here for motivation only.  It is essentially
Lemma 4 of Yu's work \cite{GY}.

Part (iii) will not be used in this paper.  However in later work we
plan to explore the application of Theorem \ref{key} to differences
between consecutive primes.  In particular we intend to use Theorem
\ref{key} to improve on Matom\"{a}ki's bound \cite{Mat}
\[\twosum{p_n\le x}{p_{n+1}-p_n\ge x^{1/2}}p_{n+1}-p_n\ll  x^{2/3}.\]
\bigskip

{\bf Acknowledgements}  This work was supported by EPSRC grant number
EP/K021132X/1. The author would like to thank the referee for pointing out
a number of misprints in an earlier version of this paper.

\section{Dirichlet Polynomials and Gaps}

In this section we will describe a procedure for bounding 
\[\#\{a_n\le x:\,a_{n+1}-a_n\ge H\}\]
from above.
It turns out that it is convenient to work on dyadic ranges $x<a_n\le
2x$, and to work with gaps of length at least $2H$, rather than
$H$.  We will assume that $H\le x^{3/4}$, and we define
\[\cl{N}(H,x)=\cl{N}=\#\{x<a_n\le 2x:\,a_{n+1}-a_n\ge 2H\}.\]
Suppose that $m$ is the smallest integer with $m>a_n/H$.  Then
$m>x/H$, and $m\le a_n/H+1\le 3x/H$. Moreover if
\[\dz:=H/(4x)\]
then 
\[(1+\dz)Hm\le \left(1+\frac{H}{4x}\right)(a_n+H)= 
a_n+H+\frac{a_n}{4x}H+\frac{H}{4x}H\]
\beql{ub}
\le a_n+H+\frac{H}{2}+\frac{H}{4}< a_n+2H\le \min(4x,a_{n+1}),
\eeq
so that
\beql{I}
[Hm,(1+\dz)Hm]\subset (a_n,a_{n+1})
\eeq
in particular. For each $a_n$ counted by $\cl{N}$ there is a
corresponding integer $m$, giving us a set $\cl{M}\subset(x/H,4x/H)$
of such integers $m$, with $\#\cl{M}=\cl{N}$.

We now define
\beql{Fdef}
F(s)=\sum_{n\le N}c_nn^{-s}=\frac{1}{k!}P(s)^{k-1}P_1(s),
\eeq
with
\[P(s):=\sum_{x^{1/k}/2<p\le x^{1/k}}p^{-s}\andd
P_1(s):=\sum_{x^{1/k}<p\le 2^{k+2}x^{1/k}}p^{-s}.\]
Here we will have $N=2^{k+2} x$.  Clearly the coefficients $c_n$ of $F(s)$ 
will be minorants for the
characteristic function of the sequence $(a_j)_1^{\infty}$
if we fix  $k>\varepsilon^{-1}$. In particular, if $n\in (x,4x]$
we will always have $c_n\le 1$, and
indeed we will have $c_n=0$ unless $n$ is a member of the sequence 
$(a_j)_1^{\infty}$. Under the above assumptions, if
$m\in\cl{M}$ then
\[\sum_{Hm<n\le (1+\dz)Hm}c_n=0,\]
since there are no elements $a_j$ in the
interval $[Hm,(1+\dz)Hm]$. We will compare the above sum with the
average over a somewhat longer interval $(Hm,(1+\dw)Hm]$, with
\[\dw=\dw(x)=\exp(-\sqrt{\log x})\ge\dz.\]

We now have
\[\sum_{Hm<n\le (1+\dw)Hm}c_n=\]
\[\frac{1}{k!}\sum_q \#\{x^{1/k}<p\le
2^{k+1} x^{1/k}:\, p\mbox{ prime, } Hm/q<p\le(1+\dw)Hm/q\},\]
where $q$ runs over products $q=p_1\ldots p_{k-1}$, counted according
to multiplicity, with $x^{1/k}/2<p_i\le x^{1/k}$. Since
$2^{1-k}x^{(k-1)/k}<q\le x^{(k-1)/k}$ and 
\[x<Hm<(1+\dw)Hm<2Hm<2(1+\dz)HM<8x,\]
by (\ref{ub}),
we see that the condition $Hm/q<p\le(1+\dw)Hm/q$ already implies that 
$x^{1/k}<p\le 2^{k+1} x^{1/k}$. Moreover the Prime Number Theorem
holds with a sufficiently good error term that we can deduce an
asymptotic formula
\[\#\{Hm/q<p\le(1+\dw)Hm/q:\, p\mbox{ prime} \}\sim \frac{\dw
  Hm}{q}\frac{k}{\log x}.\]
We then find that
\[\sum_{Hm<n\le (1+\dw)Hm}c_n\gg_k\frac{\dw x}{(\log x)^k}\]
so that
\beql{gd}
\sum_{Hm<n\le (1+\dw)Hm}c_n\ge\dw x(\log x)^{-k-1},
\eeq
say, for all $m\in(x/H,4x/H)$ and all large enough $x$.

Thus
\[\dw^{-1}\sum_{Hm<n\le (1+\dw)Hm}c_n-\dz^{-1}\sum_{Hm<n\le (1+\dz)Hm}c_n
\ge x(\log x)^{-k-1}\]
for every $m\in\cl{M}$, whence
\[\sum_{m\in\cl{M}}\left\{
\dw^{-1}\sum_{Hm<n\le (1+\dw)Hm}c_n-\dz^{-1}\sum_{Hm<n\le (1+\dz)Hm}c_n\right\}
\ge x(\log x)^{-k-1}\cl{N}.\]

We now follow the usual analysis of Perron's formula, as in Titchmarsh
\cite[Sections 3.12 and 3.19]{titch} for example. If
$0<T\le x$ we see that
\begin{eqnarray*}
\dz^{-1}\sum_{Hm<n\le (1+\dz)Hm}c_n&=&
\frac{1}{2\pi i}\int_{-iT}^{iT}F(s)\frac{(1+\dz)^s-1}{\dz s}(Hm)^sds\\
&&\hspace{2cm}\mbox{}+O\big(\dz^{-1}T^{-1}x\log x\big)
\end{eqnarray*}
and similarly 
\begin{eqnarray*}
\dw^{-1}\sum_{Hm<n\le (1+\dw)Hm}c_n&=&
\frac{1}{2\pi i}\int_{-iT}^{iT}F(s)\frac{(1+\dw)^s-1}{\dw s}(Hm)^sds\\
&&\hspace{2cm}\mbox{}+O\big(\dw^{-1}T^{-1}x\log x\big).
\end{eqnarray*}
Since $\dw\ge\dz$ we may conclude that
\begin{eqnarray*}
\lefteqn{\sum_{m\in\cl{M}}\left\{\dw^{-1}\sum_{Hm<n\le (1+\dw)Hm}c_n-
\dz^{-1}\sum_{Hm<n\le (1+\dz)Hm}c_n\right\}}\\
&=&\frac{1}{2\pi i}\int_{-iT}^{iT}F(s)\left\{\frac{(1+\dw)^s-1}{\dw s}
-\frac{(1+\dz)^s-1}{\dz s}\right\}H^sM(-s)ds\\
&&\hspace{1cm}\mbox{}+O\big(\dz^{-1}T^{-1}x(\log
x)\cl{N}\big),
\end{eqnarray*}
with
\beql{Mdef}
M(s):=\sum_{m\in\cl{M}}m^{-s}.
\eeq
We pause to remark that it is the use of this Dirichlet
polynomial $M(s)$ which is the most novel feature of the method
introduced by Yu \cite{GY}.  While the coefficients of $F(s)$ will have
some useful arithmetic structure, those of $M(s)$ do not.  None the less
it is possible to use $M(s)$ in a non-trivial way in what follows.

We now insist that $H$ satisfies
\[H\ge (\log x)^{k+3}.\]
Then the condition $T\le x$ 
will be satisfied if we choose
\beql{Tc}
T:=(\log x)^{k+3}\frac{x}{H}.
\eeq
Then
\[\dz^{-1}T^{-1}x\log x=o(x(\log x)^{-k-1}),\]
so that
\[\int_{-T}^T|F(it)M(it)|.\left|\frac{(1+\dw)^{it}-1}{\dw t}
-\frac{(1+\dz)^{it}-1}{\dz t}\right|dt
\gg x(\log x)^{-k-1}\cl{N}.\]
Here we have
\[\frac{(1+\dw)^{it}-1}{\dw t}-\frac{(1+\dz)^{it}-1}{\dz t}
\ll\min\{1\,,\,\dw (|t|+1)\},\]
so that the integral for $|t|\le\dw^{-1/8}$ contributes
\[\ll\dw^{3/4}\max_t |F(it)M(it)|\ll \dw^{3/4}x\cl{N}.\]
This is negligible compared to $x(\log x)^{-k-1}\cl{N}$ and we conclude that
\[\cl{N}\ll x^{-1}(\log x)^{k+1}\int_{\dw^{-1/8}}^T|F(it)M(it)|dt.\]

It is time to summarize our conclusions.
\begin{lemma}\label {l1}
Suppose that $(\log x)^{k+3}\le H\le x^{3/4}$.  Then
\[\cl{N}\ll x^{-1}(\log x)^{k+1}\int_{\dw^{-1/8}}^T|F(it)M(it)|dt.\]
\end{lemma}

We next show that the sum $P_1(it)$
must be relatively small in the range $\dw^{-1/8}\le t\le T$, and it is here
that it is crucial that we have removed the region
$|t|\le\dw^{-1/8}$.  We use a standard argument, see the proof of
Lemma 19 in Heath-Brown \cite{Pr}, for example.  In brief, the sum
$P_1(it)$ is
\[\frac{1}{2\pi i}\int_{\nu-it/2}^{\nu+it/2}\log\zeta(s+it)x^{s/k}
\frac{2^{(k+2)s}-1}{s}ds+O(x^{1/(2k)})+O(t^{-1}x^{1/k}\log x)\]
for $\nu=1+1/\log x$, by Perron's formula.  We can move the line of
integration to $\text{Re}(s)=1-(\log t)^{-3/4}$, say, and use the
bound
\[\log\zeta(s+it)\ll\log t,\]
which is valid in the Vinogradov--Korobov zero-free region.  Since 
\[\exp(\tfrac18\sqrt{\log x})\le t\le T\le x\] 
we deduce that
\beql{bd}
P_1(it)\ll x^{1/k}\exp\{-(\log x)^{1/9}\},
\eeq
say.
This small saving over the trivial bound will be enough for our purposes.
It shows that
\[(\log x)^{k+1}|F(it)|\ll x^{1/k}\exp\{-(\log x)^{1/10}\}|P(it)|^{k-1}\]
for $\dw^{-1/8}\le t\le T$.

We now conclude as follows.
\begin{lemma}\label{l2}
We have
\[\cl{N}\ll x^{-(k-1)/k}\exp\{-(\log x)^{1/10}\}
\int_0^T|P(it)^{k-1}M(it)|dt.\]
\end{lemma}

At this point we choose integers $a,b\ge 0$ with $a+b=k-1$ and apply
Cauchy's inequality to deduce that
\[\int_0^T|P(it)^{k-1}M(it)|dt\le\left\{\int_0^T|P(it)^b|^2dt\right\}^{1/2}
\left\{\int_0^T|P(it)^aM(it)|^2dt\right\}^{1/2}.\]
The standard mean-value
estimate for Dirichlet polynomials (Montgomery \cite[Theorem 6.1]{mont})
shows that
\[\int_0^T|P(it)|^{2b}dt\ll (T+x^{b/k})x^{b/k}.\]
We will require that $x^{b/k}\ge T$, so that the above bound is
$O(x^{2b/k})$. 
To handle the second integral we use part (iv) of Theorem \ref{key},
with $N=x^{a/k}$, which shows that
\[\int_0^T|P(it)^aM(it)|^2dt\ll x^{2a/k}R^2+R^2T+x^{a/k}RT,\]
with $R=\#\cl{M}=\cl{N}$. We now require that $x^{a/k}\ge T^{1/2}$,
in which case that above bound reduces to $O(x^{2a/k}R^2+x^{a/k}RT)$.

Comparing our estimates we now see that
\[R\ll 
x^{-(k-1)/k}\exp\{-(\log x)^{1/10}\}x^{b/k}(x^{a/k}R+x^{a/2k}R^{1/2}T^{1/2}).\]
Since $a+b=k-1$ one cannot have
\beql{fail}
R\ll x^{-(k-1)/k}\exp\{-(\log x)^{1/10}\}x^{b/k}.x^{a/k}R,
\eeq
and we conclude that
\beql{fl}
R\ll x^{-(k-1)/k}x^{b/k}.x^{a/2k}R^{1/2}T^{1/2}.
\eeq

We should stress at this point that unless (\ref{fail}) fails we can
draw no conclusion whatsoever as to the size of $R$.  It is for this
reason that our approach to Theorem \ref{t1} breaks down entirely when
$a_{n+1}-a_n\le x^{1/3}$.

The bound (\ref{fl}) leads at once to the following result.

\begin{lemma}
We have
\[\cl{N}=R\ll x^{-a/k}T\]
provided that $T^{1/2}\le x^{a/k}\le x^{(k-1)/k}T^{-1}$.
\end{lemma}

A suitable integer $a$ will exist provided that $T^{3/2}\le x^{(k-2)/k}$,
in which case we may choose $a$ so that $x^{(k-2)/k}T^{-1}\le
x^{a/k}\le x^{(k-1)/k}T^{-1}$.  This shows that
\[\cl{N}\ll x^{-(k-2)/k}T^2.\]
Recalling our choice (\ref{Tc}) we see that
\[\cl{N}\ll x^{1+2/k}H^{-2}(\log x)^{2k+6},\]
provided that $(x/H)^{3/2}\le x^{(k-3)/k}$, say. However, we chose $k$
to be an arbitrary fixed integer with $k>\ep^{-1}$. By taking $k$
suitably large we see that $\cl{N}\ll_{\ep}x^{1+\ep/2}H^{-2}$
as long as $H\ge x^{1/3+\ep}$. Clearly Theorem \ref{t1} now follows via dyadic
subdivision (for both the size of $a_n$, and the size of $a_{n+1}-a_n$).
Of course, as one decreases $\ep$ the
$x^{\ep}$-smooth numbers thin out, and the sum increases. Thus
one gets a sharper bound, for a larger quantity, on a longer range of
values $a_{n+1}-a_n$.

\section{Proof of Theorem \ref{key}}

Our proof of part (i) of Theorem \ref{key} follows ideas developed by
Montgomery \cite[Chapter 7, Theorem 3]{10}. We begin with the following
easy lemma.
\begin{lemma}\label{Lcomp}
Let $A(s)=\sum_{n\le N}\alpha_n n^{-s}$ with $\alpha_n$ real, and set
\[I_A(T)=\int_0^T|A(it)|^2dt.\]
Then
\[\frac{4}{\pi^2}I_A(T)\le T
\twosum{m,n\le N}{|\log(m/n)|\le 2\pi/T}|\alpha_m\alpha_n|,\]
and if the $\alpha_n$ are all non-negative real numbers, then
\[T\twosum{m,n\le N}{|\log(m/n)|\le 2\pi/T}\alpha_m\alpha_n\le \pi^2
I_A(T).\]
\end{lemma}

For the proof we use repeatedly the fact that the Fourier transform of
\[w(x):=\max(1-|x|\,,\,0)\]
is
\[\hat{w}(t)=\left(\frac{\sin(\pi t)}{\pi t}\right)^2.\]
Then
\[I_A(T)=\tfrac12\int_{-T}^T|A(it)|^2dt\le
\tfrac12(\tfrac{\pi}{2})^2\int_{-\infty}^{\infty}\hat{w}(t/2T)|A(it)|^2dt,\]
since $\sin(y)/y\ge 2/\pi$ for $|y|\le\pi/2$. Thus
\begin{eqnarray*}
  I_A(T)&\le &(\tfrac{\pi}{2})^2T
  \int_{-\infty}^{\infty}\hat{w}(u)|A(2Tui)|^2du\\
  &=&
(\tfrac{\pi}{2})^2 T
\twosum{m,n\le N}{|\log(m/n)|\le\pi/T}\alpha_m\overline{\alpha_n}
w(\pi^{-1}T\log(m/n))\\
&\le &
(\tfrac{\pi}{2})^2T
\twosum{m,n\le N}{|\log(m/n)|\le 2\pi/T}|\alpha_m\alpha_n|,
\end{eqnarray*}
from which the first inequality of the lemma follows.

Similarly, if the $\alpha_n$ are real and non-negative, then
\begin{eqnarray*}
\twosum{m,n\le N}{|\log(m/n)|\le 2\pi/T}
\alpha_m\alpha_n &\le &
(\tfrac{\pi}{2})^2\sum_{m,n\le N}\alpha_m\alpha_n\hat{w}
\big(\tfrac{1}{4\pi} T\log(m/n)\big)\\
&=&(\tfrac{\pi}{2})^2\int_{-\infty}^{\infty}w(u)|A(\tfrac12 Tui)|^2du\\
&=&(\tfrac{\pi}{2})^2\frac{2}{T}\int_{-\infty}^{\infty}w(2t/T)|A(it)|^2dt\\
&\le &(\tfrac{\pi}{2})^2\frac{2}{T}\int_{-T}^T|A(it)|^2dt\\
&=&\pi^2 T^{-1}I_A(T),
\end{eqnarray*}
as required.

Part (i) of Theorem \ref{key} follows at once from Lemma \ref{Lcomp}
on taking $A(s)=M(s)Q(s)$.  

We turn next to part (ii) of the theorem.
By part (i) it suffices to handle the special case in which $Q(s)$ is
\[Z(s):=\sum_{n\le N}n^{-s}.\]
\begin{lemma}
  Under the Lindel\"{o}f Hypothesis we have
  \[\cl{I}(\cl{M},Z)\ll_{\eta}N^2R^2+(NT)^{\eta}NRT,\]
 for any $\eta>0$.
\end{lemma}
By Perron's formula we have
\[Z(it)=\frac{1}{2\pi i}\int_{3/2-iN}^{3/2+iN}\zeta(s+it)N^s\frac{ds}{s}
+O(N^{1/2}).\]
Under the Lindel\"{o}f Hypothesis we can move the line of integration to
$\text{Re}(s)=\tfrac12$ to show that
\[Z(it)\ll_{\eta} N^{1/2}(N+|t|)^{\eta}+N/(1+|t|),\]
the second term coming from the pole at $s=1-it$ (if $|t|\le N$).

We now use this bound to deduce that
\[\cl{I}(\cl{M},Z)\ll_{\eta} N(N+T)^{2\eta}\int_0^T|M(it)|^2dt+\sup_t|M(it)|^2
\int_0^T N^2/(1+|t|)^2 dt.\]
By the usual mean value theorem (Montgomery
\cite[Theorem 6.1]{mont}) the first integral on the right is
$O(RT)$, whence 
\[\cl{I}(\cl{M},Z)\ll_{\eta} N(N+T)^{2\eta}RT+R^2N^2.\]
The lemma then follows on replacing $\eta$ by $\eta/2$.

The remaining parts of Theorem \ref{key} will require more effort.
We begin with the following lemma.
\begin{lemma}\label{liii1}
Let $N$ and $N_0$ be positive integers, and write
\[Z(s):=\sum_{n\le N}n^{-s}\andd Z_0(s):=\sum_{n\le N_0}n^{-s}.\]
Then
\beql{ge}
\cl{J}(\cl{M},Z) \ll N^2R^2T^{-1}+NN_0^{-1}R^2+
NN_0^{-1}\cl{J}(\cl{M},Z_0).
\eeq
\end{lemma}

In particular one sees that if $N_0=TN^{-1}$ then
\[\cl{J}(\cl{M},Z) \ll N^2R^2T^{-1}+N^2T^{-1}\cl{J}(\cl{M},Z_0).\]
This is exactly what one would expect from an examination of the
integral $\cl{I}(\cl{M},Z)$ if one were to apply the approximate
functional equation to $Z(it)$, changing its length from $N$ to
$T/N$.  However Lemma \ref{liii1} covers more general values of $N$
and $N_0$.  Moreover our proof can be adapted to handle Dirichlet
polynomials $Q(s)$
where the approximate functional equation does not apply.

For the proof we begin by writing
\[\cl{J}(\cl{M},Z)=\sum_{m_1,m_2\in\cl{M}}\#\cl{N}(m_1,m_2;N,T),\]
with
\[\cl{N}(m_1,m_2;N,T):=\{(n_1,n_2):n_1,n_2\le N,\,
|\log(m_1n_1/m_2n_2)|\le 2\pi/T\}.\]
Suppose that $m_1\ge m_2$, say.  Since $T\ge 20$ we have
\[\left|\frac{m_2n_2}{m_1n_1}-1\right|\le 8T^{-1}\]
so that
\[|n_1-n_2m_2/m_1|\le 8NT^{-1}.\]
The set 
\[\Lambda:=\{\big(n_1,n_2,T(n_1-n_2m_2/m_1)\big):\,(n_1,n_2)\in\Z^2\}\]
is a 2-dimensional lattice, with 
\[\det(\Lambda)=\sqrt{1+T^2(1+m_2^2/m_1^2)}.\]
According
to part (iii) of Heath-Brown \cite[Lemma 1]{annals} (which is based on
Lemma 4 of Davenport \cite{dav}) the lattice has a basis
$\b{e}_1$ and $\b{e}_2$ say, with
\[\text{det}(\Lambda)\ll|\b{e}_1|.|\b{e}_2|\ll\text{det}(\Lambda),\]
and such that, if $\b{e}=\lambda_1\b{e}_1+\lambda_2\b{e}_2$ then
$\lambda_i\ll |\b{e}|/|\b{e}_i|$ for $i=1,2$. Without loss of
generality we may assume that $|\b{e}_1|\le|\b{e}_2|$.

The vectors $\b{e}$ produced from points of $\cl{N}(m_1,m_2;N,T)$
will have length $|\b{e}|\le N\sqrt{66}$ so that the corresponding coefficients
$\lambda_i$ satisfy $|\lambda_i|\le c_0N/|\b{e}_i|$ for an
appropriate numerical constant $c_0$. In particular, if the pair
$m_1,m_2$ is such that $|\b{e}_2|>c_0N$, then $\lambda_2=0$.  
We will call such pairs $(m_1,m_2)$ ``bad''.
In this case
every vector in $\Lambda$ which arises from a pair $(n_1,n_2)\in
\cl{N}(m_1,m_2;N,T)$ 
 will be an integral scalar multiple of
$\b{e}_1$.  

For the remaining ``good'' pairs $(m_1,m_2)$ we have
$|\b{e}_1|,|\b{e}_2|\ll N$ so that there are $O(N/|\b{e}_i|)$
choices for $\lambda_i$.  This yields
\[\#\cl{N}(m_1,m_2;N,T)\ll \frac{N^2}{|\b{e}_1|.|\b{e}_2|}\ll N^2/T.\]
Since there are at most $R^2$ good pairs, the corresponding
contribution to (\ref{ge}) is $O(N^2R^2/T)$, which is satisfactory.

Suppose on the other hand, that $m_1,m_2$ is a bad pair, with $m_1\ge
m_2$, and that $\b{e}_1=(u_1,u_2,u_3)$ is the shorter of the two basis
vectors for $\Lambda$. If necessary we may replace $\b{e}_1$ by
$-\b{e}_1$ so that $u_1\ge 0$. 
If $\cl{N}(m_1,m_2;N,T)$ is non-empty, containing a pair
$(n_1^{(0)},n_2^{(0)})$ say, then we must have
$(n_1^{(0)},n_2^{(0)})=\lambda_0(u_1,u_2)$ for some positive
integer $\lambda_0$.  Moreover
\beql{x}
|\log(m_1n_1^{(0)}/m_2n_2^{(0)})|\le 2\pi/T,
\eeq
and hence
\[|\log(m_1u_1/m_2u_2)|\le 2\pi/T.\]
On the other hand, if this latter condition is met, then
\[(n_1,n_2)\in\cl{N}(m_1,m_2;N,T)\]
if and only if 
$(n_1,n_2)=\lambda(u_1,u_2)$ for some integer $\lambda$ satisfying
\[0<\lambda u_1\le N \andd 0<\lambda u_2\le N.\]
The condition on $\lambda$ is that it should belong to a certain
interval, $I=(0,L]$ say. 

We stress that the lattice $\Lambda$, the
basis $\b{e}_1,\b{e}_2$ and hence the interval $I$, all depend only on
$m_1$ and $m_2$, and not on $N$. However a pair $(m_1,m_2)$ may be bad
for some $N$ and good for others.

For any interval $I=(0,L]$ and any real 
$\rho>0$ one has
\[\#(\Z\cap I)\le \meas(I)=\rho^{-1}\meas(\rho I)\le
\rho^{-1}\big(1+\#(\Z\cap\rho I)\big).\]
Taking $\rho=N_0/N$ we therefore deduce that if $(m_1,m_2)$ is bad for
$N$ then
\begin{eqnarray*}
\lefteqn{\#\cl{N}(m_1,m_2;N,T)}\\
&\le& NN_0^{-1}\left(1+
\#\{\lambda\in\Z: 0<\lambda u_i\le N_0,\,(i=1,2)\}\right)\\
&\le& NN_0^{-1}\left(1+\#\cl{N}(m_1,m_2;N_0,T)\right),
\end{eqnarray*}
since every pair $(n_1,n_2)=\lambda(u_1,u_2)$ produced above will
satisfy (\ref{x}).
It follows that bad pairs $(m_1,m_2)$ contribute a total
\[\ll NN_0^{-1}R^2+NN_0^{-1}\cl{J}(\cl{M},Z_0)\] 
to (\ref{ge}), which suffices for the lemma.

We proceed to develop the above technique so as to apply to part (iv) 
of Theorem \ref{key}.  We will use the following lemma.
\begin{lemma}\label{lv1}
Let $N$ be a positive integer, and write
\[Q(s)=\frac{1}{k!}\left(\sum_{p\le N^{1/k}}p^{-s}\right)^k,\andd
Q_h(s)=\frac{1}{h!}\left(\sum_{p\le N^{1/k}}p^{-s}\right)^h,\]
  with $p$ running over primes.
Then there exists is a non-negative integer $h\le k-1$ such that
\beql{ge1}
\cl{J}(\cl{M},Q) \ll_k N^2R^2T^{-1}+ R^2+N^{(k-h)/k}\cl{J}(\cl{M},Q_h).
\eeq
\end{lemma}

For the proof it will be convenient to write
\[\mathcal{Q}_h:=\{q\in\N:\,q=p_1\ldots p_h: p_i\le N^{1/k}\}.\]
We follow the same procedure as before, writing
\[\cl{J}(\cl{M},Q)=\sum_{m_1,m_2\in\cl{M}}\#\cl{N}(m_1,m_2;N,T),\]
where now
\[\cl{N}(m_1,m_2;N,T):=\{(q_1,q_2):q_1,q_2\in\cl{Q}_k,\,
|\log(m_1q_1/m_2q_2)|\le 2\pi/T\}.\]
As previously we have
\[\#\cl{N}(m_1,m_2;N,T)\ll N^2/T\]
for good pairs $(m_1,m_2)$, contributing $O(N^2R^2/T)$ in Lemma
\ref{lv1}.

For each bad pair $(m_1,m_2)$ there will be a corresponding integer
vector $(u_1,u_2)$ such that every pair $(q_1,q_2)$ belonging to
$\cl{N}(m_1,m_2;N,T)$ takes the form
$(q_1,q_2)=\lambda(u_1,u_2)$. Clearly we must have
$u_1,u_2\in\cl{Q}_h$ and $\lambda\in\cl{Q}_{k-h}$ for some integer $h$
in the range $0\le h\le k$. We will focus attention on the value of
$h$ which makes the largest contribution. If $h=k$ then $\lambda$ must
be 1, and since there are at most $R^2$ bad pairs $(m_1,m_2)$ the
overall contribution in Lemma \ref{lv1} is $O(R^2)$. Otherwise we note
that $\#\cl{Q}_{k-h}\le N^{(k-h)/k}$, so that there
are at most $N^{(k-h)/k}$ possibilities for $\lambda$.  Moreover
\begin{eqnarray*}
\lefteqn{\#\{(m_1,m_2)\in\cl{M}^2: (u_1,u_2)\in\cl{Q}_h^2,\,
|\log(m_1u_1/m_2u_2)|\le 2\pi/T\}}\\
&\ll_h&\cl{J}(\cl{M},Q_h),\hspace{7cm}
\end{eqnarray*}
so that the overall contribution to Lemma \ref{lv1} is
$O_k(N^{(k-h)/k}\cl{J}(\cl{M},Q_h))$. This completes the proof.

\section{Completing the Proof of Theorem \ref{key}}

We have already dealt with parts (i) and (ii) of the theorem.
For part (iii) we begin by applying Lemma \ref{liii1} 
along with part (i) of Theorem \ref{key}.
These yield
\[\cl{I}(\cl{M},Z)\ll 
N^2R^2+NN_0^{-1}R^2T+NN_0^{-1}\cl{I}(\cl{M},Z_0).\]
However by Cauchy's inequality we obtain
\[\cl{I}(\cl{M},Z_0)\le\left\{\int_0^T|M(it)|^2dt\right\}^{1/2}
\left\{\int_0^T|M(it)Z_0(it)^2|^2dt\right\}^{1/2}.\]
The first integral is $O(RT)$ by the usual mean value theorem (Montgomery
\cite[Theorem 6.1]{mont}). Moreover, if 
\[D:=\max\{d(n):\,n\le N_0^2\}\ll_{\eta} N_0^{\eta},\]
then
the Dirichlet series $D^{-1}Z_0(s)^2$ has coefficients bounded by 1,
and supported in $(0,N_0^2]$.  Thus, according to part (i) of Theorem
\ref{key}, if
\[Z_1(s):=\sum_{n\le N_0^2}n^{-s},\]
then
\[\cl{I}(\cl{M},D^{-1}Z_0^2)\ll T\cl{J}(\cl{M},D^{-1}Z_0^2)\le
T\cl{J}(\cl{M},Z_1)\ll\cl{I}(\cl{M},Z_1).\]
We may therefore deduce that
\beql{tw}
\cl{I}(\cl{M},Z)\ll 
N^2R^2+NN_0^{-1}R^2T+NN_0^{-1}(RT)^{1/2}\{D^2\cl{I}(\cl{M},Z_1)\}^{1/2}.
\eeq

If we take $N_0=\sqrt{N}$ then $Z=Z_1$, and (\ref{tw}) yields
\[\cl{I}(\cl{M},Z)\ll N^2R^2+N^{1/2}R^2T+NRTD^2.\]
In particular, if $N\ge T^{2/3}$ we obtain
\beql{tt}
\cl{I}(\cl{M},Z)\ll_{\eta} N^2R^2+N^{1+\eta}RT,
\eeq
as claimed.
In the remaining case in which $N\le T^{2/3}$ we use (\ref{tw}) a 
second time, taking
$N_0=\max(T^{1/3},(RT)^{1/4})$. Then $Z_1$ has length at least
$T^{2/3}$, so that (\ref{tt}) yields
\[\cl{I}(\cl{M},Z_1)\ll_{\eta} N_0^4R^2+N_0^{2+\eta}RT.\]
Inserting this into (\ref{tw}) produces
\begin{eqnarray*}
\cl{I}(\cl{M},Z)&\ll_{\eta}&
N^2R^2+NN_0^{-1}R^2T\\
&&\hspace{1cm}\mbox{}+
NN_0^{-1}(RT)^{1/2}\{T^{2\eta}(N_0^4R^2+N_0^2RT)\}^{1/2}\\
&\ll_{\eta}&
N^2R^2+NN_0^{-1}R^2T+NR^{3/2}T^{1/2+\eta}N_0+NRT^{1+\eta}.
\end{eqnarray*}
However $NN_0^{-1}R^2T\le NR^{7/4}T^{3/4}$, and 
\[NR^{3/2}T^{1/2+\eta}N_0\le NR^{3/2}T^{5/6+\eta}+NR^{7/4}T^{3/4+\eta}.\]
For the first term on the right we have
\[NR^{3/2}T^{5/6}=(NRT)^{1/3}(NR^{7/4}T^{3/4})^{2/3}\le
\max\big(NRT\,,\,NR^{7/4}T^{3/4}\big).\]
We may therefore deduce that
\[\cl{I}(\cl{M},Z)\ll_{\eta}N^2R^2+NRT^{1+\eta}+NR^{7/4}T^{3/4+\eta}\]
when $N\le T^{2/3}$. In conjunction with (\ref{tt}) this suffices for
part (iii) of Theorem \ref{key}, since $NR^{7/4}T^{3/4+\eta}\le NRT^{1+\eta}$
when $R\le T^{1/3}$.
\bigskip

Finally we deal with part (iv) of Theorem \ref{key}. By part (i) of
the theorem (\ref{ge1}) becomes
\[\cl{I}(\cl{M},Q) \ll_k N^2R^2+ R^2T+N^{(k-h)/k}\cl{I}(\cl{M},Q_h).\]
Moreover, H\"{o}lder's inequality shows that if
\[P(s)=\sum_{p\le N^{1/k}}p^{-s}\]
then
\[\int_0^T|M(it)P(it)^h|^2dt\le
\left(\int_0^T|M(it)|^2dt\right)^{(k-h)/k}
\left(\int_0^T|M(it)P(it)^k|^2dt\right)^{h/k}.\]
As before, the first integral on the right is $O(RT)$, and we deduce
that
\[\cl{I}(\cl{M},Q) \ll_k 
N^2R^2+ R^2T+N^{(k-h)/k}(RT)^{(k-h)/k}\cl{I}(\cl{M},Q)^{h/k}.\]
Since $0\le h<k$ it then follows that
\[\cl{I}(\cl{M},Q) \ll_k N^2R^2+ R^2T+NRT,\]
as required.

\section{Theorems \ref{t2} and \ref{t3}}

Our starting point for the proof of Theorem \ref{t2} will be Lemma
\ref{l2}.  We shall assume that
\beql{st}
H\ge T^{4/k}.
\eeq
In the alternative case we have 
\[\cl{N}\ll x\ll x^{1+6/k}H^{-3/2},\]
which is satisfactory for Theorem \ref{t2} when $k$ is taken
sufficiently large.

We break the range $(0,T]$ into subintervals $(h-1,h]$, and
pick points $\beta_h, \gamma_h$ from each interval at which $|P(it)|$
and $|M(it)|$ are maximal.  We then set $b_j=\beta_{2j-1},c_j=\gamma_{2j-1}$
if the odd values of $h$ make the larger overall contribution, and
otherwise we take $b_j=\beta_{2j},c_j=\gamma_{2j}$.  Thus
\beql{in1}
\int_0^T|P(it)^{k-1}M(it)|dt\ll\sum_j|P(ib_j)|^{k-1}|M(ib_j)|,
\eeq
with points $b_j,c_j\in[0,T]$ satisfying $b_{j+1}-b_j\ge 1$ and
$c_{j+1}-c_j\ge 1$ for every relevant index $j$.

Values of $j$ for which $|P(ib_j)|\le 1$ contribute $O(T\cl{N})$ to
(\ref{in1}). We classify the remaining indices into $O(\log x)$ sets
according to the dyadic range $(V,2V]$ in which $|P(ib_j)|$
lies. Focusing on the value of $V$ which makes the largest
contribution, we re-label the relevant points as $b_j, c_j$ for $1\le
j\le J$.  We may then deduce that there is
some $V$ for which
\[\int_0^T|P(it)^{k-1}M(it)|dt\ll T\cl{N}+
(\log x)V^{k-1}\sum_{j=1}^J|M(ic_j)|,\]
with $|P(ib_j)|\ge V$ for $1\le j\le J$.  If we insert this into Lemma
\ref{l2} we see that we must have
\beql{in2}
\cl{N}\ll x^{-(k-1)/k}\exp\{-(\log x)^{1/11}\}V^{k-1}
\sum_{j=1}^J|M(ic_j)|,
\eeq
since $T\ll x^{(k-1)/k}$ under the assumption (\ref{st}).

One way to use this bound is to apply Cauchy's inequality, noting that
\[\sum_{j=1}^J|M(ic_j)|^2\ll T\cl{N}\]
by the well-known mean-value  estimate of
Montgomery \cite[Theorem 7.3]{mont} (with $Q=1, \chi=1, \delta=1$).
This yields
\[\cl{N}\ll x^{-(k-1)/k}V^{k-1}\left(JT\cl{N}\right)^{1/2},\]
and hence
\beql{in4}
\cl{N}\ll x^{-2(k-1)/k}V^{2k-2}JT.
\eeq

We proceed to estimate $J$ using the standard machinery of mean and
large values of Dirichlet polynomials. Let $M\le x$ be a parameter to
be decided, and choose an integer $\le k$ so that 
\beql{in3}
 x^{r/k}\le M<x^{(r+1)/k}.
\eeq
The Dirichlet polynomial
$A(s):=P(s)^r$ then has coefficients which are $O_k(1)$ in size, and
supported on integers up to $x^{r/k}$.  

Suppose firstly that $V\le x^{3/4k}$.  We apply Montgomery's
mean-value  estimate to $P(s)^r$, which
shows that
\[J\ll_k V^{-2r}(x^{r/k}+T)x^{r/k}.\]
Thus on taking $M=T$ we find that
\[J\ll_k V^{-2r}T^2,\]
whence the estimate (\ref{in4}) produces
\[\cl{N}\ll x^{-2(k-1)/k}V^{2k-2r}T^3.\]
Under our assumption that $V\le x^{3/4k}$ this yields
\[\cl{N}\ll x^{-2(k-1)/k}x^{3/2-3r/2k}T^3.\]
Finally, recalling that we have chosen $M=T$, we see that (\ref{in3})
produces
\[\cl{N}\ll
x^{-2(k-1)/k}x^{3/2}T^{-3/2}x^{3/2k}T^3=T^{3/2}x^{-1/2+7/2k}
\ll x^{1+4/k}H^{-3/2},\]
under the assumption (\ref{st}).
\bigskip

We turn now to the case in which $V\ge x^{3/4k}$, where we shall use
the large values estimate of Huxley, \cite[page 117]{hux} (with a
trivial modification to handle our spacing condition on the $b_k$).
In order to specify $M$ we shall define $\sigma$ by taking
$V=x^{\sigma/k}$, whence $\tfrac34\le\sigma\le 1$.  We then set 
\beql{ch}
M=\left(Tx^{2/k}\right)^{1/(4\sigma-2)}.
\eeq
In view of (\ref{Tc}) this choice will satisfy $M\le x$ provided that 
(\ref{st}) holds and $x$ is large enough.

Huxley's result now yields
\begin{eqnarray*}
J&\ll& \{V^{-2r}x^{2r/k}+V^{-6r}Tx^{4r/k}\}(\log x)^5\\
&=&\{x^{(2-2\sigma)r/k}+Tx^{(4-6\sigma)r/k}\}(\log x)^5.
\end{eqnarray*}
However $Mx^{-1/k}\le x^{r/k}\le M$ by (\ref{in3}), whence
\beql{in5}
J\ll \{M^{2-2\sigma}+TM^{4-6\sigma}x^{2/k}\}(\log x)^5\ll
M^{2-2\sigma}(\log x)^5,
\eeq
recalling our choice (\ref{ch}) for $M$.

We plan to insert this in (\ref{in2}), 
using the fact that $M(it)\ll\cl{N}$ to deduce that
\[\cl{N}\ll x^{-(k-1)/k}\exp\{-(\log x)^{1/11}\}V^{k-1}J\cl{N}.\]
We now apply the bound (\ref{in5}) together with the fact that
$V=x^{\sigma/k}$ to deduce that
\[\cl{N}\ll x^{-(1-\sigma)(k-1)/k}\exp\{-(\log x)^{1/11}\}(\log x)^5 
M^{2-2\sigma}\cl{N}.\]
Thus either $\cl{N}=0$, in which case there is nothing to prove, or 
\[\left(M^2x^{-(k-1)/k}\right)^{1-\sigma}\gg 
\exp\{(\log x)^{1/11}\}(\log x)^{-5}.\] 
In particular, if $\cl{N}\not=0$, we must have 
\[M\ge x^{(k-1)/2k}\]
so that our
definition (\ref{ch}) yields
\beql{lb}
x^{2\sigma-1}\le Tx^{(1+2\sigma)/k}\le Tx^{3/k}.
\eeq

Finally we combine the estimates (\ref{in4}) and (\ref{in5}) to deduce that
\begin{eqnarray*}
\cl{N}&\ll& T(\log x)^5x^{-2(k-1)/k}V^{2k-2}M^{2-2\sigma}\\
&=& T(\log x)^5x^{-(2-2\sigma)(k-1)/k}
\left(Tx^{2/k}\right)^{(1-\sigma)/(2\sigma-1)}\\
&=& T(\log x)^5x^{-(k-1)/(2k)}\left(Tx^{2/k}\right)^{1/2}
\left(\frac{x^{(2\sigma-1)(k-1)/k}}{Tx^{2/k}}\right)^{(4\sigma-3)/(4\sigma-2)}.
\end{eqnarray*}
According to (\ref{lb}) we have
\[\frac{x^{(2\sigma-1)(k-1)/k}}{Tx^{2/k}}\le 
\frac{x^{(2\sigma-1)}}{Tx^{2/k}}\le x^{1/k},\]
and since we are assuming that $\sigma\ge\tfrac34$ we find that
\[\cl{N}\ll T(\log
x)^5x^{-(k-1)/(2k)}\left(Tx^{2/k}\right)^{1/2}x^{1/k}
\ll T^{3/2}x^{-1/2+3/k}\ll x^{1+4/k}H^{-3/2}.\]

In every case we therefore have $\cl{N}\ll x^{1+6/k}H^{-3/2}$, and on
taking $k$ suitably large we see that Theorem \ref{t2} follows, by
dyadic subdivision of the ranges for both $a_n$ and $a_{n+1}-a_n$.
It therefore remains to establish Theorem \ref{t3}.  However we have
\begin{eqnarray*}
\twosum{a_n\le x}{a_{n+1}-a_n\le x^{1/3+\ep}}(a_{n+1}-a_n)^2
&\ll_{\ep}&  x^{1/6+\ep/2}
\twosum{a_n\le x}{a_{n+1}-a_n\le x^{1/3+\ep}}(a_{n+1}-a_n)^{3/2}\\
&\ll_{\ep}&  x^{1/6+\ep/2}\sum_{a_n\le x}(a_{n+1}-a_n)^{3/2}\\
&\ll_{\ep}&  x^{7/6+3\ep/2},
\end{eqnarray*}
by Theorem \ref{t2}. Combining this with the estimate from Theorem \ref{t1}
yields
\[\sum_{a_n\le x}(a_{n+1}-a_n)^{3/2}\ll_{\ep} x^{7/6+3\ep/2}\]
which suffices for Theorem \ref{t3}.

\bigskip

\bigskip

Mathematical Institute,

Radcliffe Observatory Quarter,

Woodstock Road,

Oxford

OX2 6GG

UK

\bigskip

{\tt rhb@maths.ox.ac.uk}

\end{document}